%% file: ms.tex
\documentclass[journal]{IEEEtran}
\usepackage[subpreambles=false]{standalone}
\usepackage{smoothing}
\usepackage{import}
\usepackage[T1]{fontenc}

\graphicspath{{img/}{sections/img/}}

\begin{document}
\title{Dynamically Iterated Filters\\
\large A unified framework for improved iterated filtering via smoothing}
	\author{Anton~Kullberg, Martin~A.~Skoglund, Isaac~Skog,~\IEEEmembership{Senior Member,~IEEE}, Gustaf~Hendeby,~\IEEEmembership{Senior Member,~IEEE}%
    \thanks{\noindent This work was partially supported by the Wallenberg AI,
			Autonomous Systems and Software Program (\textsc{WASP}) funded
			by the Knut and Alice Wallenberg Foundation.}
		\thanks{\noindent A. Kullberg, M. A. Skoglund and G. Hendeby are with the Department of Electrical Engineering, Linköping University, Linköping SE-58183, Sweden (E-mail: \{anton.kullberg, martin.skoglund, gustaf.hendeby\}@liu.se).}%
        \thanks{I. Skog is with Uppsala University, Uppsala SE-75105, Sweden (E-mail: isaac.skog@angstrom.uu.se).}
        \thanks{M. A. Skoglund is also affiliated with Eriksholm Research Center, Snekkersten, Denmark.}
	}
\IEEEpubid{}
\markboth{}{}

\maketitle

\begin{abstract}
Typical iterated filters, such as the \gls{iekf}, \gls{iukf}, and \gls{iplf}, have been developed to improve the linearization point (or density) of the likelihood linearization in the well-known \gls{ekf} and \gls{ukf}.
A shortcoming of typical iterated filters is that they do not treat the linearization of the transition model of the system.
To remedy this shortcoming, we introduce \glspl{dif}, a unified framework for iterated linearization-based nonlinear filters that deals with nonlinearities in both the transition model and the likelihood, thereby constituting a generalization of the aforementioned iterated filters.
We further establish a relationship between the general \gls{dif} and the approximate iterated Rauch-Tung-Striebel smoother.
This relationship allows for a Gauss-Newton interpretation, which in turn enables explicit step-size correction, leading to damped versions of the \glspl{dif}.
The developed algorithms, both damped and non-damped, are numerically demonstrated in three examples, showing superior mean-squared error as well as improved parameter tuning robustness as compared to the analogous standard iterated filters.
\end{abstract}


\section{Introduction}\label{sec:introduction}
\import{./sections/}{introduction}

\section{Problem Formulation}\label{sec:problemformulation}
\import{./sections/}{problem}

\section{Dynamically Iterated Filters}\label{sec:dyniter}
\import{./sections/}{dyniter}

\section{Damped Dynamically Iterated Filters}\label{sec:damped}
\import{./sections/}{damped}

\section{Numerical Examples}\label{sec:examples}
\import{./sections/}{examples}

\section{Discussion \& Conclusion}\label{sec:conclusion}
\import{./sections/}{discussion}

\appendix
\section{Identities}\label{app:identities}
\import{./sections/}{matrixidentities}
\section{Loss function}\label{app:lossfunction}
\import{./sections/}{lossequivalence}

\bibliographystyle{IEEEtran} 
\bibliography{ms}
\end{document}

%% file: sections/introduction.tex
The problem of state estimation is ubiquitous and appears in a wide range of fields, such as engineering, robotics, economics, etc.
This task requires a system model that describes the system's dynamical evolution as well as a measurement model that links observed quantities to the state of the system. 
The system model can be represented either functionally or probabilistically as a transition and measurement density. 
The two representations are given by
\begin{equation}
  \begin{aligned}
    \state_{k+1} &= \dynmod(\state_{k},\pnoise_{k}), \pnoise_k\sim p(\pnoise_k)&\tikzmarknode{space}{\qquad} p(\state_{k+1}|\state_k)\\
  \obs_k &= \obsmod(\state_{k}, \onoise_{k}),~\onoise_k\sim p(\onoise_k) &\qquad p(\obs_k|\state_k)
  \end{aligned}
\end{equation}%
\begin{tikzpicture}[overlay, remember picture]
  \draw[->] (space)++(-0.3,-0.25) to [out=0, in=-180] ++(0.3,0);
\end{tikzpicture}%
Here, $\state_k,~\obs_k,~\pnoise_k$, and $\onoise_k$ denote the state, the measurement, the process noise, and the measurement noise at time $k$, respectively.
With this model, the state estimation problem involves computing either the smooothing posterior $p(\state_{0:k}|\obs_{1:k})$ or the filtering posterior $p(\state_k|\obs_{1:k})$.
If the model is affine with additive Gaussian noise, the most well-known filtering algorithm is the analytically tractable \gls{kf}, which is the optimal estimator in the \gls{mse} sense \cite{kalmanNewApproachLinear1960}.
The smoothing analog is given by the \gls{ks}, more specifically the \gls{rts} smoother \cite{RauchRTSSmoother1965}.

\begin{figure}[tb]
  \centering
  \input{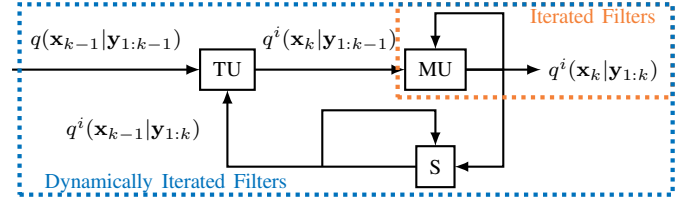}
  \caption{Schematic illustration of a \gls{dif}. 
  While standard iterated filters (dotted orange) only re-linearize the measurement update (MU), \glspl{dif} (dotted blue) also re-linearize the time update (TU) through a smoothing step (S).
  The smoothed distribution $q^i(\state_{k-1}|\obs_{1:k})$ is used to re-linearize the time update and smoothing step, whereas the current approximate posterior $q^i(\state_k|\obs_{1:k})$ is used for the measurement update.
  The steps are iterated until some convergence criterion is met.}
  \label{fig:schematic}
\end{figure}

In many practical problems, a nonlinear system model is imperative to accurately capture the system's behavior.
However, with these models, the state estimation problem becomes analytically intractable, now requiring approximate inference techniques.
Consequently, a plethora of approximate filtering algorithms have been developed, including but not limited to the \gls{ekf}, \gls{ukf}, \gls{ckf}, \gls{pf}, \gls{pmf}, and many more \cite{kalmanNewApproachLinear1960,julierNewApproachFiltering1995,arasaratnamCubatureKalmanFilters2009,schonRBPF2005,BucyPMF1971,GordonPF1993}.
Similarly, smoothing analogs such as the \gls{eks}, \gls{uks}, and \gls{ps} have also been developed \cite{CoxEKS1964,sarkkaUnscentedRauchTung2008,KitagawaParticleSmoother1996}.
In this paper, we focus on a family of methods that locally linearize the nonlinear model at each time step and subsequently utilize the \gls{kf} recursions.
This family encompasses not only the \gls{ekf} but also various sigma-point filters, as these can be interpreted as methods that employ statistical linearization \cite{lefebvreCommentNewMethod2002}.
These methods include, among others, the \gls{ukf}, the \gls{ckf}, various sigma-point filters, and the Gaussian \gls{pf} \cite{julierNewApproachFiltering1995,arasaratnamCubatureKalmanFilters2009,closasMultipleQuadratureKalman2012,steinbringLRKFRevisitedSmart2014,wangSphericalSimplexRadialCubature2014,kotechaGaussianParticleFiltering2003,wuCommentsGaussianParticle2005}.

The linearization-based methods in focus here are generally computationally efficient and are in wide practical use, with the \gls{ekf} being the {\it de~facto} standard choice in many applications.
However, linearization-based methods may suffer from convergence issues when the models are highly nonlinear.
To remedy these convergence issues, iterative versions of these methods have been developed, such as the \gls{iekf}, the \gls{iukf} and the \gls{iplf} \cite{denhamSequentialEstimationWhen1965a, Bell93,teunissenLocalConvergenceIterated1994,skoglundIterativeUnscentedKalman2019,zhanIteratedUnscentedKalman2007,garcia-fernandezIteratedStatisticalLinear2014,garcia-fernandezPosteriorLinearizationFilter2015,jazwinskiStochasticProcessesFiltering1970,sibleyIteratedSigmaPoint2006}.
These iterative methods essentially perform multiple iterations of the measurement update in every filter update, each time re-linearizing the model around the current iterate, and can thereby improve the filter convergence and by extension the posterior approximation.
Many of the methods have since been interpreted as \gls{gn} methods, whereby further improvements have been proposed.
For instance, to improve local convergence properties, \emph{damping} has been introduced, which is essentially a step-length correction based on some heuristic \cite{skoglundExtendedKalmanFilter2015,raitoharjuDampedPosteriorLinearization2018,Nocedal2006Numerical,sarkkaLevenbergMarquardtLineSearchExtended2020,lindqvistPosteriorLinearisationSmoothing2021}.
These damped methods have proven to be numerically efficient, often converge within a few iterations, and are broadly applicable to a wide range of problems.

However, previously proposed methods have focused their attention on nonlinearities in the measurement model (likelihood) and most often neglect any nonlinearities in the dynamical evolution of the system.
For some cases, this is completely reasonable, as nonlinearities in the measurement model (likelihood) affect the posterior distribution to a greater extent than the transition model (prior).
In other cases, this can lead to suboptimal performance or even divergence, as we show in our numerical experiments.
To remedy this, in \cite{kullbergIteratedFiltersNonlinear2023}, we developed a unified scheme of so-called \glsxtrfullpl{dif} that deal with the nonlinearities in both the measurement model (likelihood) \emph{as well} as in the transition model (prior).
The unified scheme encompasses both analytical and statistical linearization and thereby constitutes a generalization of the standard class of iterated linearization-based filters.
In this paper, we further elaborate on this unified scheme and extend our prior work in multiple ways.
Firstly, we provide a detailed connection between \glspl{dif} and iterated \gls{rts} smoothers, which allows us to interpret \glspl{dif} as \gls{gn} algorithms.
Secondly, in light of the \gls{gn} connection, we introduce damping to further improve the local convergence properties of the \glspl{dif}.
Lastly, we demonstrate the new family of methods on three numerical examples, each of which shows the benefits of the \glspl{dif} compared to their analogous iterated and non-iterated siblings.

The paper is organized as follows.
In \cref{sec:problemformulation}, the general state estimation problem is formulated from a probabilistic viewpoint, and linearization-based methods are particularly discussed.
\cref{sec:dyniter} derives the \glspl{dif}, which closely follows our previous work \cite{kullbergIteratedFiltersNonlinear2023}. 
In \cref{sec:damped}, the \glspl{dif} are interpreted as \gls{gn} algorithms, and damped versions thereof are developed.
\cref{sec:examples} presents multiple examples where the \glspl{dif}, as well as their damped versions, are applied to different nonlinear systems.
Lastly, \cref{sec:conclusion} discusses relevant connections to related work and outlines future research directions.

%% file: sections/schematic.tex
\begin{tikzpicture}[squarednode/.style={rectangle, draw=black, text opacity=1, thick, minimum size=5mm},
every node/.style={inner sep=5pt,outer sep=0pt, font=\footnotesize}]
\node[squarednode] (timeupdate) {TU};
\node[squarednode] (measupdate) [right=of timeupdate, xshift=1cm] {MU};
\node[squarednode] (smoothing) [below=of measupdate, yshift=.25cm] {S};
\draw[-latex, thick] ($(timeupdate.west)+(-2.5, 0)$) -- node[label=above:{$q(\state_{k-1}|\obs_{1:k-1})$}, yshift=-.1cm] {} (timeupdate.west);
\draw[-latex, thick] (timeupdate.east) -- node[label=above:{$q^i(\state_{k}|\obs_{1:k-1})$}, yshift=-.1cm] {} (measupdate.west);
\draw[-latex, thick] (measupdate.east) -| ++(0.5, 0) |- (smoothing.east);
\draw[-latex, thick] (measupdate.east) -- ++(1.0, 0) node[xshift=-.25cm, label=right:{$q^i(\state_{k}|\obs_{1:k})$}] {};
\draw[-latex, thick] ($(measupdate.east)+(0.5, 0)$) -- ++(0, .75) -| (measupdate.north);
\draw[-latex, thick] (smoothing.west) -| node[label=left:{$q^i(\state_{k-1}|\obs_{1:k})$}, yshift=.5cm] {} (timeupdate.south);
\draw[-latex, thick] ($(smoothing.west)+(-1.25, 0)$) -- ++(0, .75) -| (smoothing.north);
\draw[ultra thick, dotted, Orange] (measupdate.north west)+(-.1, .6)
node [xshift=2.4cm, yshift=-.6cm, label={[xshift=.2cm]above:{Iterated Filters}}] {}
rectangle ($(measupdate.south east)+(2.75, -0.1)$);
\draw[ultra thick, dotted, RoyalBlue] (timeupdate.north west)+(-2.4, .6)
node [xshift=1.75cm, yshift=-1.9cm, label={[xshift=.2cm]below:{Dynamically Iterated Filters}}] {}
rectangle ($(smoothing.south east)+(2.9, -0.1)$);
\end{tikzpicture}

%% file: sections/problem.tex
The general discrete-time state estimation problem is described here from a probabilistic point of view. 
We then describe linearization-based approaches in more detail and highlight their shortcomings.

\subsection{Discrete-time State Estimation}
Consider a discrete-time \gls{ssm} given by
\begin{subequations}\label{eq:ssm}
\begin{align}
\state_{k+1} &= \dynmod(\state_{k}) + \pnoise_{k}\label{eq:dyneq}, & p(\pnoise_k) &= \Ndist(\pnoise_k;\mathbf{0}, \mathbf{Q})\\
\obs_k &= \obsmod(\state_k) + \onoise_k\label{eq:obseq}, & p(\onoise_k) &= \Ndist(\onoise_k;\mathbf{0}, \mathbf{R})
\end{align}
\end{subequations}
We assume that $\state_k\in \mathcal{X},\forall k$ and that $\pnoise_k$ and $\onoise_k$ are mutually independent.
This can equivalently be expressed as 
\begin{subequations}
\begin{align}
p(\state_{k+1}|\state_{k}) &= \Ndist(\state_{k+1};\dynmod(\state_{k}), \mathbf{Q})\label{eq:transitiondensity}\\
p(\obs_k|\state_k) &= \Ndist(\obs_k;\obsmod(\state_k), \mathbf{R})\label{eq:observationdensity}.
\end{align}
\end{subequations}
Lastly, the initial state distribution is assumed to be given by
\begin{equation}
p(\state_0)=\Ndist(\state_0;\hat{\state}_{0|0}, \mathbf{P}_{0|0}).
\end{equation}

Now, given a sequence of measurements $\obs_{1:k}=\{\obs_i\}_{i=1:k}$, the smoothing problem involves computing the posterior of the state trajectory, \ie,
\begin{align}
    \label{eq:smoothingproblem}
    &p(\state_{0:k}|\obs_{1:k}) = \frac{1}{\mathbf{Z}_{1:k}} p(\state_0)\prod_{i=1}^{k} p(\obs_i|\state_i)p(\state_i|\state_{i-1})\\
    &\mathbf{Z}_{1:k}\triangleq\int_{\mathcal{X}} p(\state_0)\prod_{i=1}^{k} p(\obs_i|\state_i)p(\state_i|\state_{i-1}) ~d\state_0 \cdots d\state_k\nonumber,
\end{align}    
where $\mathbf{Z}_{1:k}$
is the marginal likelihood of $\obs_{1:k}$.
The joint smoothing distribution \cref{eq:smoothingproblem} can be analytically computed through a \gls{ks}, such as the \gls{rts} smoother, in the case of affine functions $\dynmod$ and $\obsmod$ \cite{RauchRTSSmoother1965,sarkkaBayesianFilteringSmoothing2013}.
For completeness, the computational steps of the \gls{ks} are given in \cref{alg:kalmansmoother}, where ${}_{k|K}$ indicates an estimate at time $k$ given measurements up until and including time $K$.

On the other hand, the filtering problem involves computing the marginal posteriors
\begin{align}
    \label{eq:statemarginaltimek}
&p(\state_k|\obs_{1:k}) = \frac{p(\obs_k|\state_k)\!\int_{\mathcal{X}}\! p(\state_k|\state_{k-1})p(\state_{k-1}|\obs_{1:k-1})~d\state_{k-1}}{\mathbf{Z}_{k}}\\
&\mathbf{Z}_{k} \triangleq \int_{\mathcal{X}} p(\obs_k|\state_k) p(\state_k|\state_{k-1})p(\state_{k-1}|\obs_{1:k-1})~d\state_{k-1}d\state_k,\nonumber
\end{align}
for all times $k$.
Similarly, if $\dynmod$ and $\obsmod$ are affine, the analytical solution is given by the Kalman filter \cite{kalmanNewApproachLinear1960}, which can also be found in \cref{alg:kalmansmoother}.
In the general case, neither the smoothing nor the marginal posteriors can be analytically computed.
To address this, we next discuss linearization-based approximations.


\subsection{Linearization-based Estimation}
Linearization-based estimation can generally be viewed as approximating the transition and measurement densities, \ie,
\begin{equation*}
 p(\state_{k+1}|\state_k) \overset{a}{\approx} q(\state_{k+1}|\state_k),\quad
p(\obs_k|\state_k) \overset{a}{\approx} q(\obs_k|\state_k),
\end{equation*}
where $q(\state_{k+1}|\state_k)$ and $q(\obs_k|\state_k)$ are the respective approximations.
In particular, linearization-based filters and smoothers assume affine Gaussian densities for $q(\state_{k+1}|\state_k)$ and $q(\obs_k|\state_k)$ and the \gls{kf} or \gls{ks} is then applied to this approximate model.
We will next detail the typical strategies employed in finding $q(\state_{k+1}|\state_k)$ and $q(\obs_k|\state_k)$ by means of linearization.

\subsubsection{Time Update}
When computing \cref{eq:statemarginaltimek}, two integrals need to be evaluated.
Firstly, the Chapman-Kolmogorov equation
\begin{equation}\label{eq:chapman}
p(\state_k|\obs_{1:k-1}) = \int_{\mathcal{X}} p(\state_k|\state_{k-1})p(\state_{k-1}|\obs_{1:k-1})~d\state_{k-1},
\end{equation}
which, if $p(\state_{k-1}|\obs_{1:k-1})$ is Gaussian, has a closed form solution given by \cref{eq:timeupdate}, provided that $p(\state_k|\state_{k-1})$ is Gaussian and \cref{eq:dyneq} is affine.
Given that \cref{eq:transitiondensity} is Gaussian, we thus seek an affine approximation of the transition function $\dynmod$ as
\begin{equation}\label{eq:dynapprox}
\dynmod(\state_{k-1}) \approx \mathbf{A}_\dynmod\state_{k-1} + \mathbf{b}_\dynmod + \eta_\dynmod,
\end{equation}
where $p(\eta_\dynmod) = \Ndist(\eta_\dynmod;\mathbf{0},\boldsymbol{\Omega}_\dynmod)$, capturing the error of the linearized model.
Consequently, we approximate the transition density $p(\state_{k}|\state_{k-1})$ as
\begin{equation}\label{eq:transdensityapprox}
q(\state_{k}|\state_{k-1}) = \Ndist(\state_{k}; \mathbf{A}_\dynmod\state_{k-1}+\mathbf{b}_\dynmod, \mathbf{Q}+\boldsymbol{\Omega}_\dynmod).
\end{equation}
The parameters $\mathbf{A}_\dynmod, \mathbf{b}_\dynmod$, and $\boldsymbol{\Omega}_\dynmod$ are typically chosen either through analytical or statistical linearization.
Analytical linearization yields the \gls{ekf} update, while \gls{sl} results in the various sigma-point filters, such as the \gls{ukf} with the unscented transform.

\input{\currfiledir kalmansmoother}

\subsubsection{Measurement Update}
After approximating the time update, the measurement update proceeds with computation of the (now approximate) marginal likelihood, \ie,
\begin{equation}\label{eq:normalizationapprox}
\mathbf{Z}_{k} \approx \int_{\mathcal{X}} p(\obs_k|\state_k)q(\state_k|\obs_{1:k-1})~d\state_k.
\end{equation}
Similar to \cref{eq:dynapprox}, \cref{eq:normalizationapprox} has a closed form solution if $p(\obs_k|\state_k)$ is Gaussian and \cref{eq:obseq} is affine.
Therefore, since \cref{eq:observationdensity} is Gaussian, we seek an affine approximation of the measurement function $\obsmod$ as
\begin{equation}\label{eq:obsapprox}
\obsmod(\state_k) \approx \mathbf{A}_\obsmod\state_k + \mathbf{b}_\obsmod + \eta_\obsmod,
\end{equation}
where $p(\eta_\obsmod) = \Ndist(\eta_\obsmod; \mathbf{0}, \boldsymbol{\Omega}_\obsmod)$.
Thus, the measurement density $p(\obs_k|\state_k)$ is approximated as
\begin{equation}\label{eq:obsdensityapprox}
q(\obs_k|\state_k) = \Ndist(\obs_k; \mathbf{A}_\obsmod\state_k + \mathbf{b}_\obsmod, \mathbf{R} + \boldsymbol{\Omega}_\obsmod).
\end{equation}
With \cref{eq:transdensityapprox,eq:obsdensityapprox}, the approximate marginal posterior \cref{eq:statemarginaltimek} is now given by
\begin{equation}\label{eq:statemarginaltimekapprox}
q(\state_k|\obs_{1:k}) = \frac{q(\obs_k|\state_k)q(\state_k|\obs_{1:k-1})}{\int_{\mathcal{X}} q(\obs_k|\state_k)q(\state_k|\obs_{1:k-1})~d\state_k},
\end{equation}
which is analytically tractable and given by \cref{eq:measupdate}.
Notably, analytical linearization in \cref{eq:obsapprox} around the expected value of $q(\state_{k}|\obs_{1:k-1})$ yields the \gls{ekf} measurement update, while \gls{sl} recovers the different sigma-point measurement update(s).

\subsubsection{Linearization Strategies}\label{subsec:linearization}
\input{\currfiledir linearization}
\subsection{Discussion}
The quality of the approximate marginal posterior \cref{eq:statemarginaltimekapprox} directly depends on the approximations \cref{eq:transdensityapprox,eq:obsdensityapprox}.
In turn, the quality of \cref{eq:transdensityapprox,eq:obsdensityapprox} directly depend on the 
linearization points (densities), typically chosen as the previous approximate posterior and approximate predictive distributions.
To improve the approximation of the measurement density \cref{eq:obsdensityapprox}, iterated filters, such as the \gls{iekf}, \gls{iukf} or \gls{iplf} have been proposed \cite{jazwinskiStochasticProcessesFiltering1970, garcia-fernandezIteratedStatisticalLinear2014, skoglundIterativeUnscentedKalman2019, zhanIteratedUnscentedKalman2007}.

However, except for \cite{raitoharjuPosteriorLinearisationFilter2022}, none of these algorithms address the approximate transition density \cref{eq:transdensityapprox}, despite its direct impact on the approximate marginal posterior.
Typically, this is motivated by the fact that nonlinearities in the measurement model affect the resulting posterior approximation to a greater extent than the transition model.
Nevertheless, by accounting for the linearization error in the transition model $\dynmod$, the performance of iterated filters can be improved even further.
Furthermore, standard iterated filters are not even useful in the case of a nonlinear transition function $\dynmod$ but linear measurement function $\obsmod$, even though the linearization of $\dynmod$ also affects the quality of the approximate posterior.
To address this limitation, we propose a linearization-based family of methods that encompasses both the transition density as well as the measurement density approximations.

%% file: sections/kalmansmoother.tex
\begin{algorithm}[tb]
\caption{Kalman smoother}
\label{alg:kalmansmoother}
\begin{enumerate}
\item Time update
\begin{subequations}
\label{eq:timeupdate}
\begin{align}
\hspace*{\dimexpr-\leftmargini}\hat{\state}_{k+1|k} &= \mathbf{A}_\dynmod\hat{\state}_{k|k} + \mathbf{b}_\dynmod\\
\hspace*{\dimexpr-\leftmargini}\mathbf{P}_{k+1|k} &= \mathbf{A}_\dynmod\mathbf{P}_{k|k}\mathbf{A}_\dynmod^\top + \mathbf{Q} + \boldsymbol{\Omega}_\dynmod.
\end{align}
\end{subequations}
\item Measurement update
\begin{subequations}
\label{eq:measupdate}
\begin{align}
\hspace*{\dimexpr-\leftmargini}\hat{\state}_{k|k} &= \hat{\state}_{k|k-1} + \mathbf{K}_k(\obs_k - \mathbf{A}_\obsmod\hat{\state}_{k|k-1} - \mathbf{b}_\obsmod)\\
\hspace*{\dimexpr-\leftmargini}\mathbf{P}_{k|k} &= \mathbf{P}_{k|k-1} - \mathbf{K}_k\mathbf{A}_\obsmod\mathbf{P}_{k|k-1}\\
\hspace*{\dimexpr-\leftmargini}
\mathbf{K}_k &\triangleq \mathbf{P}_{k|k-1}\mathbf{A}_\obsmod^\top(\mathbf{A}_\obsmod\mathbf{P}_{k|k-1}\mathbf{A}_\obsmod^\top + \mathbf{R} + \boldsymbol{\Omega}_\obsmod)^{-1}.
\end{align}
\end{subequations}
\item Smoothing step
\begin{subequations}
\label{eq:smoothstep}
\begin{flalign}
\hspace*{\dimexpr-\leftmargini}\hat{\state}_{k|K} &= \hat{\state}_{k|k} + \mathbf{G}_k (\hat{\state}_{k+1|K} - \hat{\state}_{k+1|k})\\
\hspace*{\dimexpr-\leftmargini}\mathbf{P}_{k|K} &= \mathbf{P}_{k|k} + \mathbf{G}_k(\mathbf{P} _{k+1|K} -\nonumber\\ &\quad~ \mathbf{A}_\dynmod\mathbf{P}_{k|k}\mathbf{A}_\dynmod^\top-\mathbf{Q}-\boldsymbol{\Omega}_\dynmod)\mathbf{G}_k^\top\\
\hspace*{\dimexpr-\leftmargini}\mathbf{G}_k &\triangleq  \mathbf{P}_{k|k}\mathbf{A}_\dynmod^\top(\mathbf{A}_\dynmod\mathbf{P}_{k|k}\mathbf{A}_\dynmod^\top + \mathbf{Q} + \boldsymbol{\Omega}_\dynmod)^{-1}
\end{flalign}
\end{subequations}
\end{enumerate}
\end{algorithm}

%% file: sections/linearization.tex
Given a nonlinear model
\begin{equation*}
\mathbf{z} = \mathbf{g}(\state),
\end{equation*}
we wish to find an affine representation
\begin{equation}
\mathbf{g}(\state)\approx \mathbf{A}\state + \mathbf{b} + \eta,
\end{equation}
where $\eta\sim\Ndist(\eta; \mathbf{0}, \boldsymbol{\Omega})$.
There are three free parameters in the affine approximation: $\mathbf{A}, \mathbf{b}$ and $\boldsymbol{\Omega}$.
Analytical linearization, through first-order Taylor expansion, selects the parameters as
\begin{equation}\label{eq:analyticallinearization}
\mathbf{A}=\frac{d}{d\state}\mathbf{g}(\state)|_{\state=\bar{\state}},
\quad \mathbf{b} = \mathbf{g}(\state)|_{\state=\bar{\state}} - \mathbf{A}\bar{\state},
\quad \boldsymbol{\Omega}=\mathbf{0},
\end{equation}
where $\bar{\state}$ is the point about which the function $\mathbf{g}(\state)$ is linearized.
Note that $\boldsymbol{\Omega}=\mathbf{0}$ essentially implies that the linearization is assumed to be error-free, even though there exist higher-order alternatives.

\Glsxtrfull{sl} instead linearizes w.r.t. a distribution $p(\state)$, which in practice can be viewed as linear regression with samples from $p(\state)$ transformed through $\mathbf{g}(\state)$.
Assuming that such a distribution $p(\state)=\Ndist(\state;\hat{\state},\mathbf{P})$ is given, \gls{sl} selects the affine parameters as
\begin{subequations}\label{eq:statisticallinearization}
\begin{align}
\mathbf{A} &= \Psi^\top \mathbf{P}^{-1}\\
\mathbf{b} &= \bar{\mathbf{z}}-\mathbf{A}\hat{\state}\\
\boldsymbol{\Omega} &= \Phi-\mathbf{A} \mathbf{P} \mathbf{A}^\top\\
\bar{\mathbf{z}} &= \mathbb{E}[\mathbf{g}(\state)]\\
\Psi &= \mathbb{E}[(\state-\hat{\state})(\mathbf{g}(\state) - \bar{\mathbf{z}})^\top]\\
\Phi &= \mathbb{E}[(\mathbf{g}(\state) - \bar{\mathbf{z}})(\mathbf{g}(\state) - \bar{\mathbf{z}})^\top],
\end{align}
\end{subequations}
where the expectations are taken w.r.t. $p(\state)$.
The major difference from analytical linearization is that typically ${\boldsymbol{\Omega}\neq 0}$, implying that the linearization error is captured.
Typically, the expectations in \cref{eq:statisticallinearization} are not analytically tractable and practically, one often resorts to some numerical integration technique, such as the unscented transform.
Note that the parameter choices in \cref{eq:statisticallinearization} can be derived by minimizing the \gls{mse} between the function $\mathbf{g}(\state)$ and the corresponding affine approximation.
A more comprehensive discussion is out of the scope of this paper, but for additional details see, \eg, \cite{elishakoffStochasticLinearization2017}.

%% file: sections/dyniter.tex
We here recapitulate the main points of our derivation in \cite{kullbergIteratedFiltersNonlinear2023}, where we first derived the \glspl{dif}.
To that end, note that to optimally approximate both the transition density \cref{eq:transdensityapprox}, as well as the measurement density \cref{eq:obsdensityapprox} at time $k$, we naturally need to seek an approximate posterior over both $\state_{k-1}$ and $\state_k$.


\subsection{Optimal Affine Approximations}
Let two auxiliary variables $\auxil_k$ and $\auxil_{k-1}$ be defined as
\begin{subequations}
\begin{align}
\auxil_{k-1} &= \dynmod(\state_{k-1}) + \psi, & p(\psi) &= \Ndist(\mathbf{0}, \alpha\mathbf{I})\\
\auxil_k &= \obsmod(\state_k) + \phi, & p(\phi) &= \Ndist(\mathbf{0}, \beta\mathbf{I})
\end{align}
\end{subequations}
where $\psi$ and $\phi$ are independent of each other as well as the process noise $\pnoise_k$ and the measurement noise $\onoise_k$.
Note that $\auxil_{k-1}\to \dynmod(\state_{k-1})$ and $\auxil_k\to\obsmod(\state_k)$ as $\alpha,~\beta\to 0$.
The true joint posterior of $\state_{k-1},~\state_k,~\auxil_{k-1}$ and $\auxil_k$ is then given by
\begin{align}\label{eq:jointdecomposition}
&\fulljoint \nonumber \\ \propto{}&\fulljointexpand.
\end{align}
We assume that the approximate posterior decomposes in a similar manner, \ie,
\begin{align}\label{eq:approximatejoint}
 &\approxjoint \nonumber\\ \approx{}&\approxjointexpand.
\end{align}
Here, $\theta$ are the parameters of the affine approximation of the transition model and measurement model, \ie, ${\theta=[\mathbf{A}_\dynmod,\mathbf{b}_\dynmod,\boldsymbol{\Omega}_\dynmod, \mathbf{A}_\obsmod,\mathbf{b}_\obsmod,\boldsymbol{\Omega}_\obsmod]}$.

The parameters $\theta$ should then be chosen such that $\approxjoint$ closely approximates $\fulljoint$.
Specifically, the optimal parameters $\theta^*$, and thus the optimal affine approximations of $\dynmod$ and $\obsmod$, are found by minimizing a loss $\Loss$, \ie,
\begin{align}\label{eq:fullopt}
\theta^* = & \argmin_\theta\Loss.
\end{align}
The loss $\Loss$ is chosen as the \gls{kl} divergence between the true and approximate posterior, \ie
\begin{multline}
\!\!\!\!\!\Loss\triangleq\mathrm{KL}\bigl(\fulljoint \Vert \approxjoint \bigr)\\
=\mathrm{KL}\bigl( \statejoint \Vert \approxstatejoint \bigr) \\
+\mathbb{E}\bigl[\mathrm{KL}\bigl(\auxilcond{k}\Vert\approxauxilcond{k}\bigr)\bigr] \\
+\mathbb{E}\bigl[
\mathrm{KL}\bigl(\auxilcond{k-1}\Vert\approxauxilcond{k-1}\bigr)
\bigr], \label{eq:fulloss}
\end{multline}
where the expectations in $\Loss$ are taken w.r.t. the true joint posterior $\statejoint$.
See \cite{kullbergIteratedFiltersNonlinear2023} for the derivation of \cref{eq:fulloss}.
Note that $\Loss$ can be decomposed into three distinct terms, each dealing with a specific factor of \cref{eq:approximatejoint}.
The first term ensures that the true and approximate joint posterior of the state at time $k$ and $k-1$ are close.
The remaining two terms minimize the discrepancy between the affine approximations and the true measurement and transition model, respectively.

Since the expectations in $\Loss$ are taken w.r.t. the true joint posterior $\statejoint$, minimizing the loss \cref{eq:fulloss} is, unfortunately, impractical.
Nevertheless, an iterative approach can be used to approximately solve this minimization problem.

\subsection{Iterative Solution}
In order to practically optimize \cref{eq:fullopt}, we assume access to an $i$\textsuperscript{th} approximation to the state joint posterior $\statejoint\approx\approxstatejoint[i]$.
Subsequently, we substitute $\approxstatejoint[i]$ for $\statejoint$ in \cref{eq:fulloss} and thus optimize an approximate loss, given by
\begin{multline*}
    \Loss \approx \mathcal{L}(\theta^i) = \mathrm{KL}\left( \approxstatejoint[i] \Vert \approxstatejoint[i+1] \right) \\
    +\mathbb{E}_{\approxstatejoint[i]}\left[\mathrm{KL}(\auxilcond{k}\Vert\approxauxilcond[i+1]{k})\right] \\
    +\mathbb{E}_{\approxstatejoint[i]}
    \left[
    \mathrm{KL}(\auxilcond{k-1}\Vert\approxauxilcond[i+1]{k-1})
    \right],
\end{multline*}
where the expectations are now over $\approxstatejoint[i]$.
Sufficiently close to a fixed point, the first \gls{kl} term is approximately 0, simplifying the approximate loss to 
\begin{multline}
    \mathcal{L}(\theta^i) \approx \mathbb{E}_{\approxstatejoint[i]}\biggl[\mathrm{KL}\bigl(\auxilcond{k}\Vert\approxauxilcond[i+1]{k}\bigr) \\
    +\mathrm{KL}\bigl(\auxilcond{k-1}\Vert\approxauxilcond[i+1]{k-1}\bigr)
    \biggr].
\end{multline}
Technically, the optimal $\theta^*$ is given by \gls{sl} of $\dynmod$ and $\obsmod$ w.r.t. the current approximation $\approxstatejoint[i]$, see \eg, \cite{garcia-fernandezPosteriorLinearizationFilter2015}.
This only requires the marginal $\approxstate{i}{k-1}$ for \gls{sl} of $\dynmod$, and similarly, the marginal $\approxstate{i}{k}$ for \gls{sl} of $\obsmod$.
Thus, conceptually, the solution amounts to iterating a prediction, update, and smoothing step in order to provide new linearization points (densities) for both the transition density as well as the measurement density simultaneously.
These steps are repeated until fixed-point convergence, ultimately yielding an approximate posterior $q(\state_{k-1:k}|\obs_{1:k})$.
The overall algorithm is schematically depicted in \cref{fig:schematic} and is given in \cref{alg:dyniterfilter}.
\input{\currfiledir base_algorithm}

\subsection{Illustration}
\begin{figure}[tb]
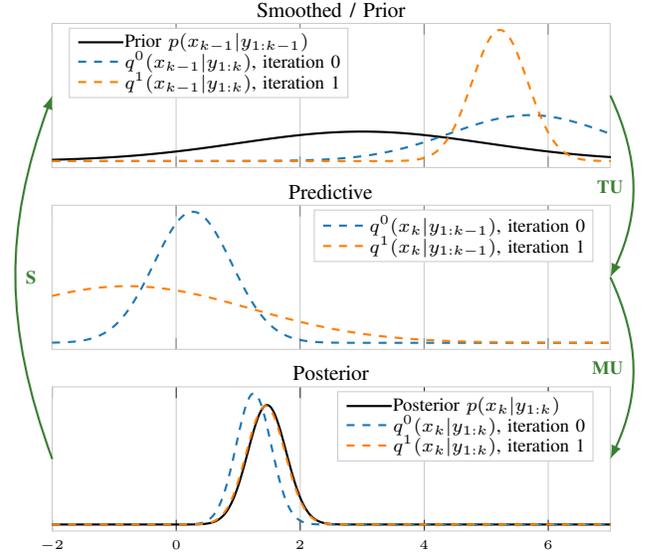

\centering
\includestandalone[]{\currfiledir img/dyniter_illustration_new}
\caption{Illustration of the (extended) \gls{dif}.
The filter moves from the prior (top) to the predictive (middle), via the time update (TU), to the posterior (bottom), via the measurement update (MU), and back up to the smoothed (top), via the smoothing step (S).
Notice that iteration 0 exactly corresponds to an \gls{ekf}.}
\label{fig:illustrativeexample}
\end{figure}

We visually illustrate the \glspl{dif} by applying the analytically linearized \gls{dif}, the \gls{diekf}, to a one-dimensional model given by
\begin{align*}
\state_{k+1} &\sim \Ndist(\state_{k+1}; a\state_k^3, Q)\\
\obs_k &\sim \Ndist(\obs_k; \state_k, R)\\
p(\state_{k-1}|\obs_{1:k-1}) &= \Ndist(\state_{k-1}; 3, 4),
\end{align*}
where $a=0.01, Q=0.1$ and $R=0.1$.
The approximate densities produced by the \gls{diekf} are visualized in \cref{fig:illustrativeexample}, where two iterations are enough to yield a good posterior approximation.
The true posterior is found simply by evaluating the posterior density over a dense grid.

\subsection{Discussion}\label{subsec:dyniterdiscussion}
Unlike regular iterated filters, the \gls{dif} proves valuable across all potential combinations of models with linear and nonlinear functions, $\dynmod$ and $\obsmod$.
It is worth noting that the \gls{dif} functions more or less like an \gls{iplf}, but also encompasses the transition density.
However, an \q{extended} version, \ie, the \gls{diekf}, is recovered by opting for analytical, rather than statistical, linearization.
Further, similar to the \gls{iukf} in \cite{skoglundIterativeUnscentedKalman2019}, an unscented version is recovered by maintaining fixed covariance matrices $\mathbf{P}_{k-1|k}^i=\mathbf{P}_{k-1|k-1}$ and $\mathbf{P}_{k|k}^i=\mathbf{P}^i_{k|k-1}$ and updating these only during the last iteration of \cref{alg:dif}.
Moreover, the \glspl{dif} may also be interpreted as \q{local} iterated smoothers, analogous to the \gls{ieks} \cite{bellIteratedKalmanSmoother1994a} and the \gls{ipls} \cite{garcia-fernandezIteratedPosteriorLinearization2017}, operating on just one time-instance and measurement.
Thus, the method can be viewed as an iterated one-step fixed-lag smoother as well, a connection which we shall now firmly establish.



%% file: sections/base_algorithm.tex
\begin{figure}[t]
\begin{minipage}{\columnwidth}
\begin{algorithm}[H]
  \caption{Dynamically iterated filter}
  \label{alg:dyniterfilter}
  \begin{algorithmic}[1]
  \Require $q(\state_{k-1}|\obs_{1:k-1})=\Ndist(\state_{k-1};\hat{\state}_{k-1|k-1},\mathbf{P}_{k-1|k-1})$
  \State Linearize $\dynmod$ about $\hat{\state}_{k-1|k-1}, \mathbf{P}_{k-1|k-1}$ by \cref{eq:analyticallinearization} or \cref{eq:statisticallinearization}
  \State Compute $\hat{\state}_{k|k-1}^0, \mathbf{P}_{k|k-1}^0$ by \cref{eq:timeupdate}
  \State Linearize $\obsmod$ about $\hat{\state}_{k|k-1}^0, \mathbf{P}_{k|k-1}^0$ by \cref{eq:analyticallinearization} or \cref{eq:statisticallinearization}
  \State Compute $\hat{\state}_{k|k}^0, \mathbf{P}_{k|k}^0$ by \cref{eq:measupdate}
  \State Compute $\hat{\state}_{k-1|k}^0, \mathbf{P}_{k-1|k}^0$ by \cref{eq:smoothstep}\vspace{.5em}
  \State $\hat{\state}_{k-1|k}, \mathbf{P}_{k-1|k}, \hat{\state}_{k|k}, \mathbf{P}_{k|k} \leftarrow $ \cref{alg:dif} or \cref{alg:dampeddif}\label{alg:dif:algchoice}
  \State \textbf{return} $\hat{\state}_{k|k}, \mathbf{P}_{k|k}$ (and $\hat{\state}_{k-1|k}, \mathbf{P}_{k-1|k}$)
  \end{algorithmic}
\end{algorithm}

\vspace*{-\baselineskip}
\begin{algorithm}[H]
    \caption{Inner loop of \gls{dif}}
\label{alg:dif}
\begin{algorithmic}[1]
    \Require $\hat{\state}_{k-1|k-1},\mathbf{P}_{k-1|k-1},\hat{\state}^0_{k-1|k},\mathbf{P}^0_{k-1|k},\hat{\state}_{k|k}^0,\mathbf{P}^0_{k|k}$
    \State $i\gets 0$
    \While{not converged}
        \State $\bar{\state}_{k-1}, \bar{\mathbf{P}}_{k-1}, \bar{\state}_{k}, \bar{\mathbf{P}}_{k} \leftarrow \hat{\state}_{k-1|k}^{i}, \mathbf{P}_{k-1|k}^{i}, \hat{\state}_{k|k}^{i}, \mathbf{P}_{k|k}^{i}$
        \State $\hat{\state}_{k-1|k}^{i+1}, \mathbf{P}_{k-1|k}^{i+1}, \hat{\state}_{k|k}^{i+1}, \mathbf{P}_{k|k}^{i+1} \leftarrow$ \cref{alg:singlestepdif}
        \State $i\gets i+1$
    \EndWhile
\end{algorithmic}
\end{algorithm}

\vspace*{-\baselineskip}
\begin{algorithm}[H]
\caption{Single step of a \gls{dif}}
\label{alg:singlestepdif}
\begin{algorithmic}[1]
    \Require $\hat{\state}_{k-1|k-1}, \mathbf{P}_{k-1|k-1}, \bar{\state}_{k-1}, \bar{\mathbf{P}}_{k-1}$, $\bar{\state}_{k}, \bar{\mathbf{P}}_{k}$
    \State Linearize $\dynmod$ about $\bar{\state}_{k-1}, \bar{\mathbf{P}}_{k-1}$ by \cref{eq:analyticallinearization} or \cref{eq:statisticallinearization}
    \State Compute $\hat{\state}_{k|k-1}, \mathbf{P}_{k|k-1}$ by \cref{eq:timeupdate}
    \State Linearize $\obsmod$ about $\bar{\state}_{k}, \bar{\mathbf{P}}_{k}$ by \cref{eq:analyticallinearization} or \cref{eq:statisticallinearization}
    \State Compute $\hat{\state}_{k|k}, \mathbf{P}_{k|k}$ by \cref{eq:measupdate}
    \State Compute $\hat{\state}_{k-1|k}, \mathbf{P}_{k-1|k}$ by \cref{eq:smoothstep}
    \State \textbf{return} $\hat{\state}_{k-1|k}, \mathbf{P}_{k-1|k}, \hat{\state}_{k|k}, \mathbf{P}_{k|k}$
\end{algorithmic}
\end{algorithm}
\end{minipage}
\end{figure}

%% file: sections/damped.tex
For more complex, highly nonlinear, models, the \glspl{dif} may sometimes diverge, similar to standard iterated filters, such as the \gls{iekf} and \gls{iukf}.
To lessen these issues, damped versions of iterated filters have been developed, where an explicit step--size has been introduced into the iterate updates, improving the (local) convergence properties \cite{skoglundIterativeUnscentedKalman2019,teunissenLocalConvergenceIterated1994,raitoharjuDampedPosteriorLinearization2018}.
Inspired by the successes in the development of damped iterated filters, we shall now apply a similar reasoning to develop damped \glspl{dif}, to reap similar benefits in terms of convergence properties.
To that end, we first establish a connection between \glspl{dif} and iterated approximate \gls{rts} smoothers.
This will allow us to interpret the \glspl{dif} as \gls{gn} methods, which in turn lets us introduce explicit step--sizes in the iterate updates, thereby improving local convergence properties of the \glspl{dif}.

\subsection{DIF as a GN Algorithm}
The loss function \cref{eq:fulloss} is suitable to motivate the iterative procedure of the \glspl{dif}.
However, as the expected \gls{kl}--terms in \cref{eq:fulloss} are not analytically tractable, it is not wise to use \cref{eq:fulloss} to select a step--size.
Thus, we require another criteria for step--size selection and establish the following proposition.

\begin{proposition}[Loss function]\label{prop:loss}
In each time step $k$, \glspl{dif} minimize
\begin{multline}
2\mathcal{L}(\state_{k-1},\state_k) = \lVert \state_{k-1} - \state_{k-1|k-1} \rVert^2_{\bP_{k-1|k-1}} \\
+ \lVert \obs_k - \obsmod(\state_k) \rVert^2_{\bR+\bOmega_\obsmod} + \lVert \state_k - \dynmod(\state_{k-1}) \rVert^2_{\bQ+\bOmega_\dynmod},\label{eq:lossfunction}
\end{multline}
where $\lVert \mathbf{v} \rVert_\mathbf{S}^2 = \mathbf{v}^\top \mathbf{S}^{-1} \mathbf{v}$.
\end{proposition}
\begin{IEEEproof}
    See \cref{app:lossfunction}.
\end{IEEEproof}

\cref{prop:loss} establishes the loss function that \glspl{dif} minimize, a dynamical analog to the loss functions of the \gls{iekf}, \gls{iukf} and \gls{iplf}.
There is, however, a caveat. 
Notice that \cref{eq:lossfunction} depends on $\bOmega_\dynmod,\bOmega_\obsmod$, \ie, the linearization error of $\dynmod$ and $\obsmod$.
In the analytically linearized \gls{diekf}, this is not an issue as the linearization error is not explicitly dealt with.
However, in the \gls{diukf} and \gls{diplf}, the linearization errors are recomputed in each iteration, meaning that the weighting of the loss changes
\begin{multline}
2\mathcal{L}^i(\state_{k-1}, \state_k) = \lVert \state_{k-1} - \state_{k-1|k-1} \rVert^2_{\bP_{k-1|k-1}} \\
+ \lVert \obs_k - \obsmod(\state_k) \rVert^2_{\bR+\bOmega_\obsmod^i} + \lVert \state_k - \dynmod(\state_{k-1}) \rVert^2_{\bQ+\bOmega_\dynmod^i},\label{eq:lossfunction_omega}
\end{multline}
which does not necessarily yield a monotonic cost decrease for decreasing residuals (note the dependence of $\bOmega$ on $i$).
For the \gls{diukf}, this is not a problem, as the state covariance is already \q{frozen} until the last iterate update.
For the \gls{diplf} on the other hand, we propose to remedy this problem by constructing an inner and outer loop, similar but not identical to that of \cite{raitoharjuDampedPosteriorLinearization2018}.
The inner loop is identical to that of the \gls{diukf}, where only the estimated mean is updated.
The outer loop on the other hand updates the state covariance and re--runs the inner loop.
This is then repeated until some convergence criteria is fulfilled, \eg, the \gls{kl} divergence between two successive approximations is negligible.
Hence, the damped \gls{diplf} is essentially an iterated damped \gls{diukf}, \ie, a \q{doubly} iterated \gls{diukf}.
The general damped \gls{dif} is given in \cref{alg:dyniterfilter}, with the choice of \cref{alg:dampeddif} on \cref{alg:dif:algchoice}.
In particular, \cref{alg:dampeddif} is the version \emph{without} the outer loop.
The only difference is that \cref{alg:dampeddif} is rerun several times, with new covariances $\mathbf{P}_{k-1|k}^0$ and $\mathbf{P}_{k|k}^0$ each time, given by the last iteration.

\input{\currfiledir damped_algorithm.tex}

\subsection{Step-size correction}
Through \cref{prop:loss}, we can establish that the \glspl{dif} are closely related to approximate one-step fixed-lag \gls{rts} smoothers.
As such, we can interpret the \glspl{dif} as \gls{gn} methods and can therefore introduce step-size correction to improve local convergence properties.

To that end, let the \gls{gn} iterate be defined by $\state^i=[\state_{k-1}^{i,\top},\state_k^{i,\top}]^\top$.
A step-size, $\alpha^i$, can then be introduced to the iterate update as 
$$
\state^{i+1} = \state^i +\alpha^i\mathbf{p}^i,
$$ 
with the step $\mathbf{p}^i$. 
The step is given by
\begin{equation}\label{eq:gnstep}
    \mathbf{p}^i = -({\mathbf{J}^i}^\top \mathbf{J}^i)^{-1}{\mathbf{J}^i}^\top \mathbf{r}^i,
\end{equation}
where 
\begin{equation}
    \mathbf{J}^i = \frac{d \mathbf{r}^i(\mathbf{s})}{d \mathbf{s}}\bigg|_{\mathbf{s}=\state^i}.
\end{equation}
Here, 
\begin{equation}
    \mathbf{r}(\state^i)=\begin{bmatrix}
        \bP_{k-1|k-1}^{-1/2}(\state_{k-1}-\state_{k-1|k-1})\\
        (\bR+\bOmega_\obsmod^\star)^{-1/2}(\obs_k-\obsmod(\state_{k}))\\
    (\bQ+\bOmega_\dynmod^\star)^{-1/2}(\state_{k}-\dynmod(\state_{k-1}))
    \end{bmatrix},
\end{equation}
where $\star$ corresponds to the two different $\bOmega_\obsmod$ and $\bOmega_\dynmod$ in \cref{eq:lossfunction}, and \cref{eq:lossfunction_omega}, respectively. 
The step $\mathbf{p}^i$ can equivalently be defined as
\begin{equation*}
    \mathbf{p}^i = \xi^{i+1} - \state^i,
\end{equation*}
which we use in \cref{alg:dampeddif}.
Here, the expression for $\xi^{i+1}$ is given immediately by a single step of the \gls{dif}, \ie, \cref{alg:singlestepdif}.
This can be seen as proposing a new iterate $\xi^{i+1}$ through the \gls{ks} in \cref{alg:kalmansmoother} and weighing its influence by $\alpha^i$.

Due to, \eg, noise, poor Hessian approximation, or strong non-linearities, a full step $\alpha^i=1$, \ie, choosing $\xi^{i+1}$ as the new iterate, may not yield a cost decrease, \ie,
$$
\mathcal{L}^\star(\state^{i}+\alpha^i \mathbf{p}^i)\nleq \mathcal{L}^\star(\state^i).
$$
Given the cost functions and control over the step-size one can find an $\alpha^i\in (0,1]$ such that 
$$
\mathcal{L}^\star(\state^i+\alpha^i\mathbf{p}^i)\leq \mathcal{L}^\star(\state^i),
$$ 
which is called a line-search. 
For more elaborate line-search strategies, we refer to, \eg, \cite{Nocedal2006Numerical}.

\subsection{Discussion}
Our construction for the \gls{diplf} is heavily inspired by \cite{raitoharjuDampedPosteriorLinearization2018}.
There is, however, a distinct difference in our construction of the inner loop.
Particularly, \cite{raitoharjuDampedPosteriorLinearization2018} freeze $\bOmega$ in their inner loop construction, in order to guarantee a monotonic cost decrease within the inner loop.
In contrast, we argue that this is not conceptually sound and can in fact hinder convergence.
Within the inner loop, the model is still re-linearized in each step and the linearization error $\bOmega$ is used in the iterate update.
Therefore, freezing it does not take the new linearization point (distribution) into account.
Further, a connection was recently established between iterated statistical linearization and quasi-Newton \cite{kullbergRelationshipsPosteriorQuasi2023}.
In particular, in \cite{kullbergRelationshipsPosteriorQuasi2023} it was shown that the \gls{iplf}, \gls{iukf} and \gls{ickf} are in fact quasi-Newton methods with very particular choices of Hessian correction.
Therefore, excluding the linearization error $\bOmega$ from the inner loop may actually lead to suboptimal step lengths.
However, we should point out that from our practical experience, this change to the inner loop seems to make little difference in practice.
There may exist problems where this does make a difference, but we are yet to find one. 
Our changes to the inner loop are thus purely for conceptual reasons and we shall not investigate it further in this paper.



%% file: sections/damped_algorithm.tex
\begin{algorithm}[t]
\begin{algorithmic}
    \caption{Inner loop of the damped \gls{diukf}/\gls{diekf}}
    \label{alg:dampeddif}
    \Require $\hat{\state}_{k-1|k-1},\mathbf{P}_{k-1|k-1},\hat{\state}^0_{k-1|k},\mathbf{P}^0_{k-1|k},\hat{\state}_{k|k}^0,\mathbf{P}^0_{k|k}$
    \State $i \leftarrow 0$
    \While{not converged}
        \State $\bar{\state}_{k-1}, \bar{\mathbf{P}}_{k-1},\bar{\state}_{k}, \bar{\mathbf{P}}_{k} \leftarrow \hat{\state}_{k-1|k}^{i}, \mathbf{P}_{k-1|k}^{0},\hat{\state}_{k|k}^{i}, \mathbf{P}_{k|k}^{0}$
        \State $\xi_{k-1}, -, \xi_{k}, - \leftarrow$ \cref{alg:singlestepdif}\vspace{.5em}\Comment{Propose new iterates}
        \State $\state^i = \begin{bmatrix} \hat{\state}^i_{k-1|k} \\ \hat{\state}^i_{k|k} \end{bmatrix},\quad \mathbf{p}^i = \begin{bmatrix}
            \xi_{k-1} - \hat{\state}^i_{k-1|k}\\ 
            \xi_{k} - \hat{\state}^i_{k|k}
        \end{bmatrix}, \quad \alpha \leftarrow 1$\vspace{.5em}
        \While{$\mathcal{L}(\state^i+\alpha \mathbf{p}^i) \leq \mathcal{L}(\state^i)$}
            \State $\alpha \leftarrow \alpha/\eta$ \Comment{Other strategies possible}
        \EndWhile\vspace{.5em}
        \State Set $\begin{bmatrix}
            \hat{\state}^{i+1}_{k-1|k}\\
            \hat{\state}^{i+1}_{k|k}
        \end{bmatrix} \leftarrow \state^i + \alpha \mathbf{p}^i$\vspace{.5em}
        \State $i \leftarrow i+1$
    \EndWhile
    \State $\bar{\state}_{k-1}, \bar{\mathbf{P}}_{k-1},\bar{\state}_{k}, \bar{\mathbf{P}}_{k} \leftarrow \hat{\state}_{k-1|k}^{i}, \mathbf{P}_{k-1|k}^{0},\hat{\state}_{k|k}^{i}, \mathbf{P}_{k|k}^{0}$
    \State $-, \mathbf{P}^{i}_{k-1|k}, -, \mathbf{P}^{i}_{k|k} \leftarrow $ \cref{alg:singlestepdif} \Comment{Covariance update}
  \State \textbf{return} $\bar{\state}_{k-1}, \bar{\mathbf{P}}_{k-1}, \bar{\state}_{k}, \bar{\mathbf{P}}_{k}$
  \end{algorithmic}
\end{algorithm}
  


%% file: sections/examples.tex
Next, we present three numerical examples illustrating the application of \glspl{dif} and their damped counterparts.
The first example consists of an artificial highly nonlinear one-dimensional problem that motivates the usage of damping in the \glspl{dif}.
The second example is a simulated maneuvering target tracking problem where the dynamics are nonlinear but the measurement model is \emph{linear}.
The third example is an acoustic localization problem where both dynamics and measurements are nonlinear.
Note that, in the first and third examples, we focus solely on the \q{extended} variants of the \glspl{dif} to keep the numerical results uncluttered.
In the second example, we consider the \gls{diekf}, \gls{diukf}, and \gls{diplf}, but not their respective damped versions, again to reduce clutter.

\subsection{One-Dimensional Nonlinear}
\begin{figure}[t]
\centering
\includegraphics[width=\columnwidth]{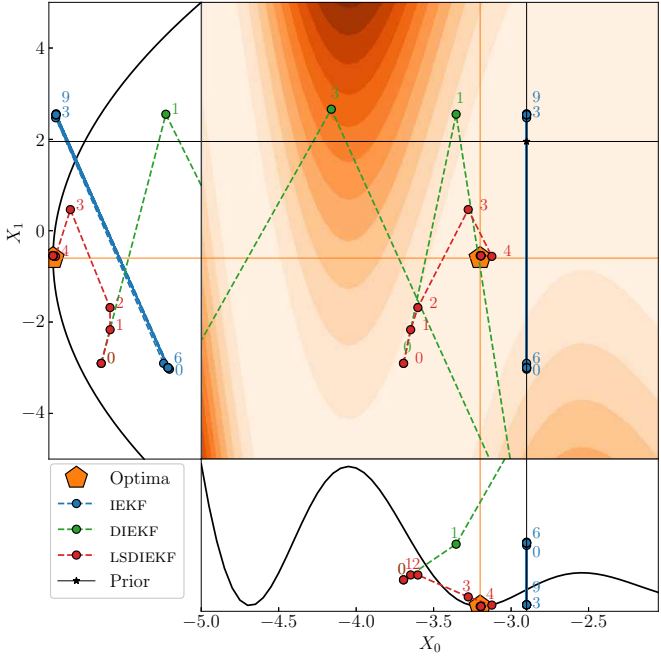}
\caption{The loss landscape of the first example as well as the iterates produced by the \gls{iekf}, \gls{diekf} and \gls{lsdiekf}. The contour plot illustrates the 2D loss landscape over the states $X_0$ and $X_1$, where a lighter orange corresponds to a lower loss. Further, the two one-dimensional subplots illustrate the marginal loss curves at the optimal $X_0$ and $X_1$, respectively.}
\label{fig:dampedmotivation}
\end{figure}

In this example, we consider an artificial system given by
\begin{subequations}
    \begin{align}
        \state_{k+1} &= \cos(\state_k)\sin(\state_k)\state_k^2 + \pnoise_k\label{eq:artificialdynamics}\\
        \obs_k &= \arctan(\state_k) + \onoise_k\label{eq:artificialobservations},
    \end{align}
\end{subequations}
with $\pnoise_k\sim\mathcal{N}(0, 0.1)$ and $\onoise_k\sim\mathcal{N}(0, 1)$.
The true state and the prior at time $0$ are given by $x=-3.2$ and $x_{0|0}\sim\mathcal{N}(-2.9, 1)$, respectively. 
We run the \gls{iekf}, \gls{diekf} and a \gls{lsdiekf} for 10 iterations each.
We visualize the 2D loss landscape in \cref{fig:dampedmotivation} over the state at time $0$, $X_0$, and time $1$, $X_1$, and overlay the iterates produced by the \gls{iekf}, \gls{diekf} and \gls{lsdiekf}, respectively.
The two one-dimensional subplots illustrate the marginal loss curves at the optima, \ie, the bottom subplot corresponds to the loss over $X_0$ at $X_1=X_1^*$ where $X_1^*$ is the optima at time $1$.
The black star and lines correspond to the prior means and the orange star and lines correspond to the optima.
From \cref{fig:dampedmotivation}, it is clear that the \gls{diekf} (green), diverges completely.
Further, the ordinary \gls{iekf} (blue), ends up in a sort of \q{limit cycle} behaviour, jumping between two modes.
On the other hand, the \gls{lsdiekf} (red), converges to the optima.

\subsection{Maneuvering Target Tracking}
For completeness, we here consider the same numerical example as in \cite{kullbergIteratedFiltersNonlinear2023}, and therefore recapitulate only the most essential information of the example.
We refer to \cite{kullbergIteratedFiltersNonlinear2023} for detailed information.


We consider a target maneuvering in a plane and describe the target state using the state vector $\state_k=\begin{bmatrix}p^x_k & v^x_k & p^y_k & v^y_k & \omega_k \end{bmatrix}^\top$. 
Here, $p^x_k,~p^y_k,~v^x_k$, and $v^y_k$ are the Cartesian coordinates and velocities of the target, respectively. Further, $\omega_k$ is the turn rate.
The transition model is given by
\begin{equation}\label{eq:maneuveringmodel}
    \state_{k+1} = \dynmod(\state_k) + \pnoise_k,
\end{equation}
where
\begin{equation*}
    \dynmod(\state_k) = \begin{bmatrix}
        1 & \frac{\sin(T \omega_k)}{\omega_k} & 0 & -\frac{(1-\cos(T \omega_k ))}{\omega_k} & 0\\
        0 & \cos(T \omega_k ) & 0 & -\sin(T \omega_k ) & 0\\
        0 & \frac{(1-\cos(T \omega_k ))}{\omega_k} & 1 & \frac{\sin(T \omega_k )}{\omega_k} & 0\\
        0 & \sin(T \omega_k ) & 0 & \cos(T \omega_k ) & 0\\
        0 & 0 & 0 & 0 & 1
        \end{bmatrix}\state_k,
\end{equation*}
and $T$ is the sampling period. 
Further, $\pnoise_k \sim \Ndist(\pnoise_k;\mathbf{0},\mathbf{Q})$ is the process noise at time $k$, with
\begin{equation*}
\mathbf{Q} =
\mathrm{blkdiag}\left( 
\begin{bmatrix} q_1 \frac{T^3}{3} & q_1\frac{T^2}{2}\\
q_1\frac{T^2}{2} & q_1T \end{bmatrix},  
\begin{bmatrix} q_1 \frac{T^3}{3} & q_1\frac{T^2}{2}\\
q_1\frac{T^2}{2} & q_1T \end{bmatrix}, q_2 \right),
\end{equation*}
where $q_1$ and $q_2$ are tunable parameters of the model.
We assume a linear measurement model of the Cartesian position of the target with additive Gaussian noise given by $\onoise_k\sim\Ndist(\onoise_k; \mathbf{0}, \sigma^2\mathbf{I})$. 

\begin{figure*}[t]
    \centering
    \includegraphics[width=\textwidth]{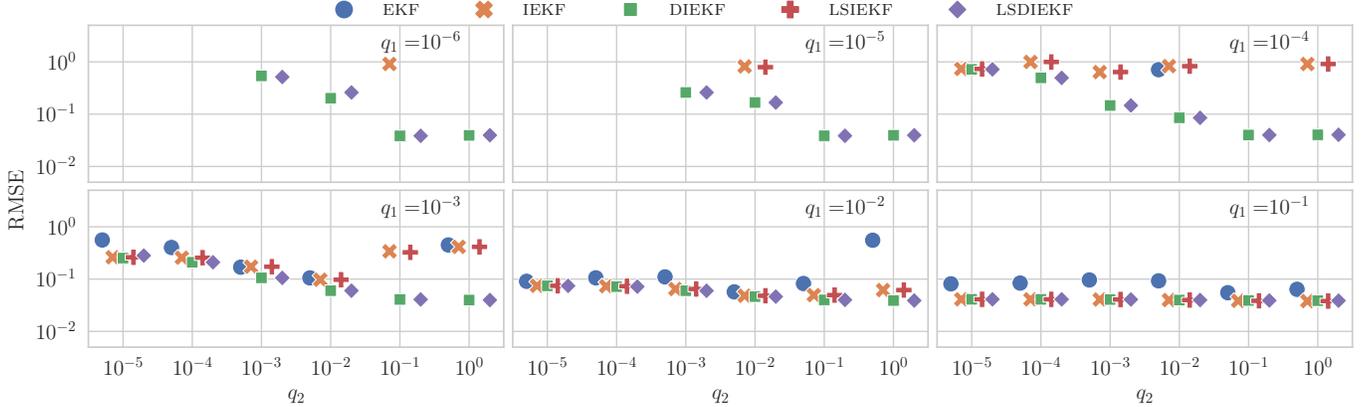}
    \caption{Position \gls{rmse} for the different methods and noise settings. Each subplot corresponds to a different value of $q_1$, indicated by the text in each subplot. An \gls{rmse} higher than approximately $\SI{1}{\meter}$ corresponds to a \q{divergent} filter based on visual inspection of resulting estimate trajectories and is left out of the plots. Note that the markers are slightly spread out to aid visibility.}
    \label{fig:tdoarmse}
\end{figure*}

\begin{figure*}[tb]
\centering
\subfloat[\label{fig:positionrmse}Position \gls{rmse}. 
The analytically linearized \gls{diekf} improves the most over its baseline the \gls{ekf} which diverges for 22/25 configurations. 
The \gls{diukf} and \gls{diplf} have similar performance, where the \gls{diplf} is slightly better for some configurations.]{
\begin{tikzpicture}
\node[label={[yshift=-.75cm, xshift=.4cm]Position \gls{rmse} $\left[\SI{}{\meter}\right]$}, label=below:{Relative \gls{rmse} Iterated/Baseline}] {\includegraphics[width=.97\textwidth]{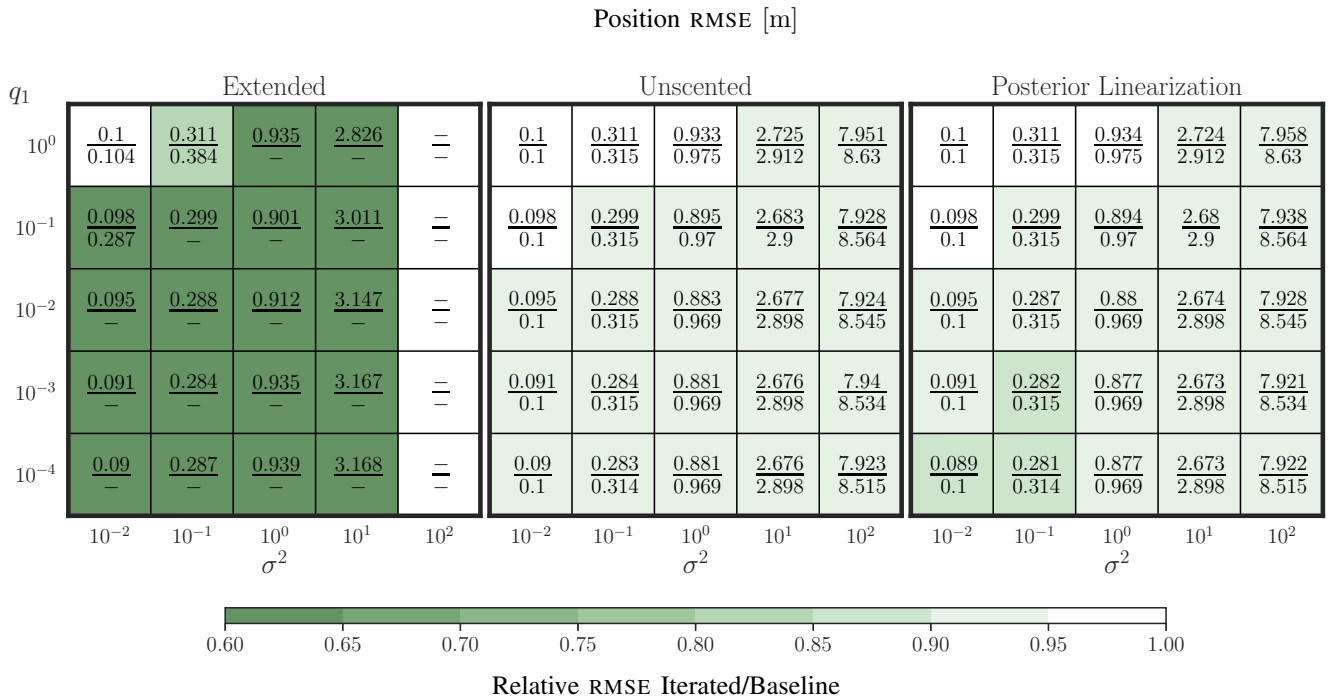}};
\end{tikzpicture}
}\\%
\subfloat[\label{fig:velocityrmse} Velocity \gls{rmse}. 
All three \glspl{dif} have substantially better performance than their corresponding baselines. 
The \gls{diekf} has 10-fold improvements and the \gls{diukf} and \gls{diplf} have approximately a 5-fold improvement in low noise regimes.
The results are even better for the \gls{diukf} and \gls{diplf} in high noise regimes, with approximately 10-fold improvements.]{
\begin{tikzpicture}
\node[label={[yshift=-.75cm, xshift=.35cm] Velocity \gls{rmse} $\left[\SI{}{\meter\per\second^{}}\right]$}, label=below:{Relative \gls{rmse} Iterated/Baseline}] {\includegraphics[width=.97\textwidth]{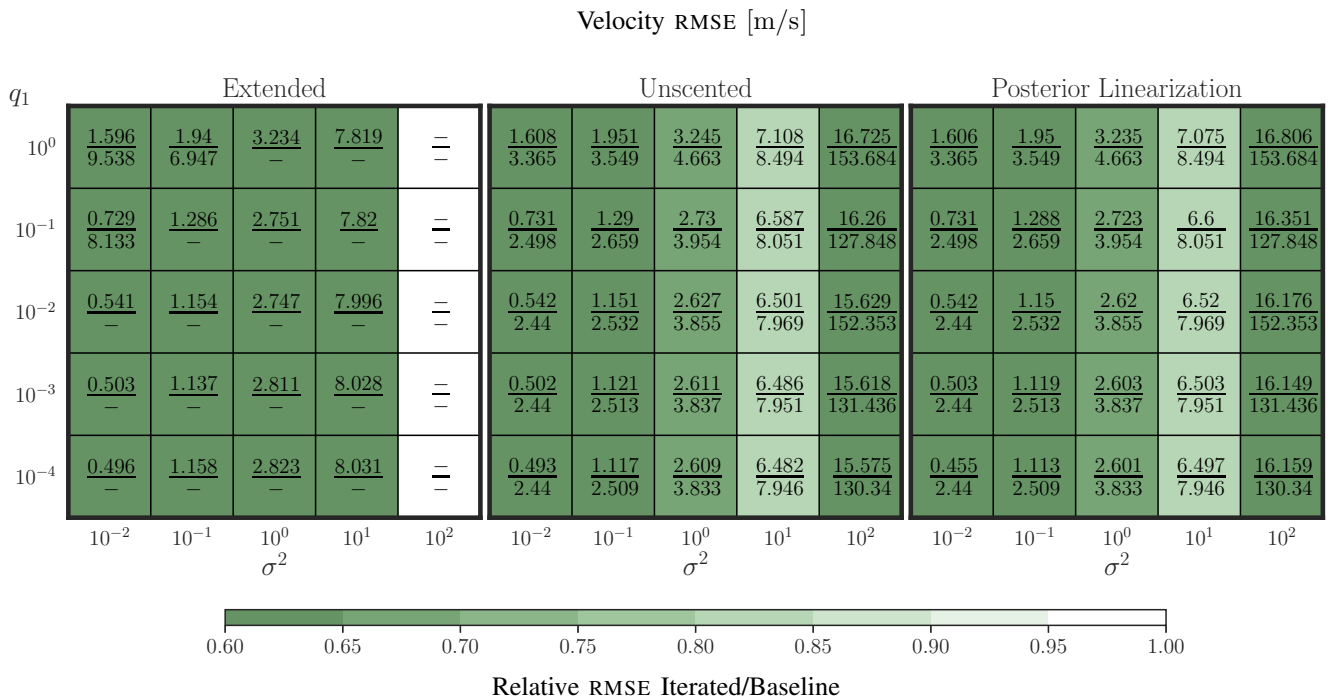}};
\end{tikzpicture}
}
\caption{\gls{rmse} for the \glspl{dif} as compared to their respective baselines in the second numerical example.
The top number in each cell is the \gls{rmse} for the \gls{dif} whereas the bottom number is the \gls{rmse} for the corresponding baseline.
A \q{$-$} indicates that the position \gls{rmse} of the filter is larger than $\sigma$ and it is thus considered to have diverged.
The left plot in \protect\subref*{fig:positionrmse} and \protect\subref*{fig:velocityrmse} shows the \gls{rmse} for the \gls{diekf} and \gls{ekf}.
The middle figure shows the \gls{diukf} and \gls{ukf} and finally, the rightmost figure shows the \gls{diplf} and \gls{ukf}.
Each cell is colored according to the \gls{rmse} of the \gls{dif} relative to the baseline; the relative \gls{rmse} is simply the fraction of the \gls{rmse}.}
\label{fig:rmse}
\end{figure*}

We consider $25$ different noise configurations, by fixing $q_2=\num{e-2}$ and varying $q_1$ and $\sigma^2$.
We then evaluate the average \gls{rmse}, separately for position and velocity, over $200$ simulations for each of the filters and their corresponding baselines.
%
%
The results are presented as $5\times 5$ matrices, where each cell corresponds to a particular noise configuration for a particular pair of algorithms, \eg, the results for the \gls{diekf} and \gls{ekf} are summarized in one matrix.
The position and velocity \glspl{rmse} are presented in \cref{fig:positionrmse} and \cref{fig:velocityrmse}, respectively.

All of the \glspl{dif} improve upon their baselines, but the analytically linearized \gls{diekf} improves the most, see \cref{fig:positionrmse}.
Surprisingly, the \gls{ekf} diverges in $22/25$ configurations whereas the \gls{diekf} lowers that to $5/25$, only diverging in the high noise scenario ($\sigma^2=\num{e+2}$).
The \gls{diukf} and \gls{diplf} improve the most in low process noise regimes, with otherwise relatively modest improvements.
However, all of the \glspl{dif} substantially improve in velocity \gls{rmse}, see \cref{fig:velocityrmse}.
For low process noise regimes the improvement is up to 10-fold for the \gls{diekf} and 5-fold for the \gls{diukf} and \gls{diplf}.
Even for modest noise levels, the \gls{diukf} and \gls{diplf} roughly manage a 2-fold performance improvement.
For the high noise scenario ($\sigma^2=\num{e+2}$), the \gls{diukf} and \gls{diplf} show a 10-fold performance improvement and bring the velocity \gls{rmse} down to reasonable levels where the \gls{rmse} for the \gls{ukf} is very high.
It is important to note that even though the measurement model is \emph{linear}, we still see significant performance improvements as well as robustness improvements for the \gls{diekf} in particular.

\subsection{Acoustic Localization}\label{subsec:acousticlocalization}
We consider an acoustic localization problem where the system has both nonlinear dynamics as well as a nonlinear measurement model.
In order to not clutter the results, we only consider the \gls{ekf}, \gls{iekf}, \gls{diekf} as well as the damped \gls{iekf} and \gls{diekf}.

An RC-car is driving along a pre-defined trajectory in an enclosed area, with dynamics modeled by \cref{eq:maneuveringmodel}.
Whilst driving, the car emits a known sound pulse at equidistant points in time.
Four microphones are used to construct \gls{tdoa} measurements.
Each measurement at microphone $j$ is thus modeled as
\begin{equation}\label{eq:tdoa}
    \obs_k^j = r^1_k - r^j_k + \onoise_k^j,\quad j=1,\dots,3
\end{equation}
where $r^j_k\triangleq \lVert \begin{bmatrix}p^x_k & p^y_k\end{bmatrix}^\top - s^j \rVert$, and $s^j$ is the 2D position of each microphone, which is assumed known.
Note that \cref{eq:tdoa} is parameterized in distance and not in time.

Before data collection, a calibration experiment was conducted to construct an error model for each microphone, and thus, $\onoise_k = \begin{bmatrix}\onoise_k^1 & \onoise_k^2 & \onoise_k^3\end{bmatrix}^\top\sim\Ndist(\mathbf{0}, \mathbf{R})$ where 
\begin{equation*}
    \mathbf{R} = \sigma_1^2\mathbf{1}\mathbf{1}^\top + \mathrm{diag}\left( \sigma_2^2~\sigma_3^2~\sigma_4^2 \right),
\end{equation*}
where $\sigma_j^2$ is the variance of microphone $j$.

We set $q_1=10^{-j},~q_2=10^{-l}$ and let $j=-6,\dots,0,~l=-5,\dots,0$, and sweep over all such pairs, \ie, $42$ different process noise configurations.
For each configuration, we compute the \gls{rmse} against a ground truth trajectory, obtained from a high--precision IR-marker positioning system.

The positional \gls{rmse} per noise configuration is presented in \cref{fig:tdoarmse}.
The \gls{ekf} diverges for almost all process noise configurations (divergence approximately corresponds to an \gls{rmse} higher than $\SI{1}{\meter}$).
The \gls{iekf} and damped \gls{iekf} improve the situation somewhat but still diverge in a lot of cases.
The \gls{diekf} and damped \gls{diekf} improve the performance even further and essentially only diverge for very low process noise tuning.
As the process noise is increased, the difference between the algorithms decreases but still favors the \gls{diekf} and damped \gls{diekf}, showing strong robustness to the process noise tuning.

%% file: sections/discussion.tex
Finally, we will discuss relationships to other works and conclude the article with some future research directions.

\subsection{Related work}
In this work, we have presented a unified framework of so-called \glsxtrfullpl{dif}.
However, certain special cases of \glspl{dif} have been presented previously in the literature.
For instance, the \gls{diplf} first presented in \cite{raitoharjuPosteriorLinearisationFilter2022}, for the purpose of robustly estimating the process noise variance, a motivation which is completely different from ours.
A similar motivation to ours appears in \cite{garcia-fernandezIteratedPosteriorLinearization2017} where the L-scan \gls{iplf} is introduced. The L-scan \gls{iplf} also re-linearizes the dynamical model through an $L$-lag smoothing procedure.
However, even in the case of $L=1$, it requires access to previous measurements, something we do not consider here, as this technically makes it a smoother.

The origin of the \gls{diekf} can be traced back to \cite{wishnerEstimationStateNoisy1968}, where it was developed through principles of feedback, lacking a proper theoretical motivation.
In contrast, we unify \glspl{dif} based on analytical and statistical linearization in a common framework. This allows for a plethora of special variants of \glspl{dif} to be developed, depending on the chosen variant of \gls{sl}.
As such, while we only consider the unscented transform as our choice of \gls{sl}, dynamically iterated analogs of other variants of sigma--point filters are readily available by making another choice of \gls{sl}.
One may also note that while we do recognize the vast literature on sigma--point filters, we have limited our references to a select few representative examples.

As alluded to in \cref{subsec:dyniterdiscussion}, the \glspl{dif} may also be interpreted as a one-step fixed-lag smoother.
As such, the different special cases have strong connections to the \gls{ieks} \cite{bellIteratedKalmanSmoother1994a} and the \gls{ipls} \cite{garcia-fernandezIteratedPosteriorLinearization2017}.
Further, the damped versions have strong connections to fixed-lag alternatives of the damped versions of the \gls{ieks} and \gls{ipls} \cite{sarkkaLevenbergMarquardtLineSearchExtended2020,lindqvistPosteriorLinearisationSmoothing2021}.

Other alternatives to nonlinear filtering are also noteworthy, such as the \gls{pf} \cite{GordonPF1993} and \gls{pmf} \cite{BucyPMF1971}, which are both applicable to a wider range of models than the algorithms presented here.
However, these have exponential complexity in the dimension of the state and are thus typically limited to lower-dimensional systems.
If the model is restricted to contain a linear Gaussian substructure, the exponential complexity has been somewhat remedied by Rao--Blackwellized versions, such as the \gls{rbpf} \cite{schonRBPF2005,KarlssonComplexityRBPF2005}.
Nevertheless, with a high-dimensional nonlinear evolution, the complexity remains.

\subsection{Summary \& Conclusion}
\input{\currfiledir conclusions}

%% file: sections/conclusions.tex
A unified framework of \glspl{dif} has been presented.
Contrary to previous iterated filters, \glspl{dif} are useful for all possible combinations of (Gaussian) linear and nonlinear transition and measurement models.
Further, the \glspl{dif} were shown to be \gls{gn} methods and damped versions of the \glspl{dif} could thus be readily developed.
The \glspl{dif} were evaluated against both their non-iterated as well as standard iterated counterparts.
In all numerical examples, the \glspl{dif} showed superior performance and robustness properties.
As a future research direction, it would be interesting to study the convergence properties of the \glspl{dif} as compared to regular iterated filters which should be possible through the lens of \gls{gn} optimization. 
Further, quantifying the properties of the models that benefit from the \gls{dif} as opposed to standard filters could be interesting.

%% file: sections/matrixidentities.tex
\begin{identity}[Woodbury Matrix Identity]\label{id:woodbury}
Assume that $\mathbf{L}$ and $\boldsymbol{\Sigma}$ are invertible matrices. Then,
\begin{align*}
&\left( \mathbf{A}\mathbf{L}\mathbf{A}^\top + \boldsymbol{\Sigma} \right)^{-1} = \\
&\boldsymbol{\Sigma}^{-1} - \boldsymbol{\Sigma}^{-1}\mathbf{A}\left( \mathbf{L}^{-1} + \mathbf{A}^\top\boldsymbol{\Sigma}^{-1}\mathbf{A} \right)^{-1}\mathbf{A}^\top\boldsymbol{\Sigma}^{-1}.
\end{align*}
\end{identity}

\begin{identity}[Kalman gain]\label{id:kalmangain}
Assume that $\mathbf{L}$ and $\boldsymbol{\Sigma}$ are invertible matrices. Then, by \cref{id:woodbury}
\begin{align*}
&\mathbf{L}\mathbf{A}^\top\left( \mathbf{A}\mathbf{L}\mathbf{A}^\top + \boldsymbol{\Sigma} \right)^{-1} = \\
&\mathbf{L}\mathbf{A}^\top\left( \boldsymbol{\Sigma}^{-1} - \boldsymbol{\Sigma}^{-1}\mathbf{A}\left( \mathbf{L}^{-1} + \mathbf{A}^\top\boldsymbol{\Sigma}^{-1}\mathbf{A} \right)^{-1}\mathbf{A}^\top\boldsymbol{\Sigma}^{-1} \right) = \\
&\overbrace{\mathbf{L}\left( \left[ \mathbf{L}^{-1}+\mathbf{A}^\top\boldsymbol{\Sigma}^{-1}\mathbf{A} \right] - \mathbf{A}^\top\boldsymbol{\Sigma}^{-1}\mathbf{A} \right)}^{=\mathbf{I}} \times \\
&\left[ \mathbf{L}^{-1}+\mathbf{A}^\top\boldsymbol{\Sigma}^{-1}\mathbf{A} \right]^{-1}\mathbf{A}^\top\boldsymbol{\Sigma}^{-1} = \\
&\left[ \mathbf{L}^{-1}+\mathbf{A}^\top\boldsymbol{\Sigma}^{-1}\mathbf{A} \right]^{-1}\mathbf{A}^\top\boldsymbol{\Sigma}^{-1}
\end{align*}
\end{identity}

%% file: sections/lossequivalence.tex
Assume a transition model and a measurement model given by
\begin{subequations}
\begin{align}
\state_{k+1} &= \dynmod(\state_k) + \pnoise_k\label{eq:transitionmodel}\\
\obs_k &= \obsmod(\state_k) + \onoise_k\label{eq:observationmodel},
\end{align}
\end{subequations}
where $\pnoise_k\sim\Ndist(\pnoise_k; \mathbf{0}, \bQ)$ and $\onoise_k\sim\Ndist(\onoise_k; \mathbf{0}, \bR)$ are mutually independent.

The (linearized) mean measurement update and mean smoothing step for this model are given by
\begin{subequations}
\begin{align}
\state_{k|k} &= \state_{k|k-1} + \bK_k \left( \obs_k - \obsmod(\state_{k|k-1}) \right)\label{eq:measurementupdate}\\
\state_{k|k+1} &= \state_{k|k} + \bG_k
\left( \state_{k+1|k+1} - \dynmod(\state_{k|k}) \right)\label{eq:smoothingstep},
\end{align}
\end{subequations}
where
\begin{align*}
\bH \triangleq& \frac{d}{d\state} \obsmod(\state)\vert_{\state=\state_{k|k-1}} & \btR &\triangleq \bR+\bOmega_\obsmod\\
\bF \triangleq& \frac{d}{d\state} \dynmod(\state)\vert_{\state=\state_{k|k}} & \btQ &\triangleq \bQ+\bOmega_\obsmod,
\end{align*}
where $\bOmega_\dynmod$ and $\bOmega_\obsmod$ are parameters of the linearization, as given by \cref{eq:statisticallinearization} or \cref{eq:analyticallinearization}.
Further, the gains $\bK_k$ and $\bG_k$ are given by
\begin{align}
\bK_k \triangleq& \bP_{k|k-1}\bH^\top\left( \bH\bP_{k|k-1}\bH^\top + \btR \right)^{-1}\nonumber\\
&\overbrace{\left( \bP_{k|k-1}^{-1} + \bH^\top\btR^{-1}\bH \right)^{-1}}^{\triangleq \Gamma}\bH^\top\btR^{-1}\label{eq:kalmangain}\\
\bG_k \triangleq& \bP_{k|k}\bF^\top\left( \bF\bP_{k|k}\bF^\top + \btQ \right)^{-1}\nonumber\\
&\left( \bP_{k|k}^{-1} + \bF^\top\btQ^{-1}\bF \right)^{-1}\bF^\top\btQ^{-1} = \nonumber\\
&\underbrace{\left( \bP_{k|k-1}^{-1} + \bH^\top\btR^{-1}\bH + \bF^\top\btQ^{-1}\bF \right)^{-1}}_{\triangleq \Lambda }\bF^\top\btQ^{-1}\label{eq:smoothinggain},
\end{align}
by \cref{id:kalmangain}.

Note that \cref{eq:measurementupdate,eq:smoothingstep} can thus also be written as
\begin{subequations}
\begin{align}
\state_{k|k} =& \state_{k|k-1} + \Gamma\bH^\top\btR^{-1} \left( \obs_k - \obsmod(\state_{k|k-1}) \right)\label{eq:measurementupdateinverted}\\
\state_{k|k+1} =
&\state_{k|k} + \Lambda\bF^\top\btQ^{-1}%
\left( \state_{k+1|k+1} - \dynmod(\state_{k|k}) \right)\label{eq:smoothingstepinverted}.
\end{align}
\end{subequations}

Using \cref{eq:measurementupdateinverted} in \cref{eq:smoothingstepinverted} yields
\begin{align}
\state_{k|k+1} = \state_{k|k-1} + %
&\overbrace{\Gamma\bH^\top\btR^{-1}%
\left( \obs_k - \obsmod(\state_{k|k-1})\right)}^{\term{term:A}} + \nonumber\\%
&\underbrace{\Lambda\bF^\top\btQ^{-1}\left( \state_{k+1|k+1} - \dynmod(\state_{k|k}) \right)}_{\term{term:B}} \label{eq:combinedsteps}
\end{align}

Focusing on \cref{term:B}, note that
\begin{align*}
&\state_{k+1|k+1} - \dynmod(\state_{k|k}) = \\%
&(\state_{k+1|k+1} - \dynmod(\state_{k|k-1})) + %
(\dynmod(\state_{k|k-1}) - \dynmod(\state_{k|k})) \approx \\%
&(\state_{k+1|k+1} - \dynmod(\state_{k|k-1})) + ~\\%
&(\tikzmarknode[strike out,draw]{1}{\dynmod(\state_{k|k})} + \bF(\state_{k|k-1} - \state_{k|k}) - \tikzmarknode[strike out,draw]{1}{\dynmod(\state_{k|k})}) = \\
&(\state_{k+1|k+1} - \dynmod(\state_{k|k-1})) %
- \bF\Gamma\bH^\top\btR^{-1}\left( \obs_k-\obsmod(\state_{k|k-1}) \right),
\end{align*}
where the last equality comes from \cref{eq:measurementupdateinverted}.

Thus, \cref{term:B} can be written as
\begin{align}
&\underbrace{\Lambda\bF^\top\btQ^{-1}\left( \state_{k+1|k+1} - \dynmod(\state_{k|k-1})\right)}_{\term{term:C}} +~\nonumber\\%
&\underbrace{- \Lambda\bF^\top\btQ^{-1}%
\bF\Gamma\bH^\top\btR^{-1}\left( \obs_k-\obsmod(\state_{k|k-1}) \right)}_{\term{term:D}}.
\end{align}

Combining \cref{term:A,term:D} yields
\begin{align}\label{eq:obspart}
\Lambda \bF^\top\btQ^{-1}\bF%
&\overbrace{\left[ \left( \Lambda \bF^\top\btQ^{-1}\bF \right)^{-1} \Gamma \bH^\top\btR^{-1} - \Gamma\bH^\top\btR^{-1}\right]}^{\term{term:E}}\times\nonumber\\
&\left( \obs_k-\obsmod(\state_{k|k-1}) \right).
\end{align}

We proceed with \cref{term:E} which can be simplified as
\begin{align*}
&\left( \Lambda \bF^\top\btQ^{-1}\bF \right)^{-1} \Gamma \bH^\top\btR^{-1} - \Gamma\bH^\top\btR^{-1} = \\%
&\left( \bF^\top\btQ^{-1}\bF \right)^{-1} \left( \bP_{k|k-1}^{-1} + \bH^\top\btR^{-1}\bH + \tikzmarknode[strike out,draw]{1}{\bF^\top\btQ^{-1}\bF} \right) \Gamma \bH^\top\btR^{-1} - \\
&~\tikzmarknode[strike out,draw]{1}{\Gamma\bH^\top\btR^{-1}} = \\%
&\left( \bF^\top\btQ^{-1}\bF \right)^{-1} \overbrace{\left( \bP_{k|k-1}^{-1} + \bH^\top\btR^{-1}\bH \right)}^{=\Gamma^{-1}} \Gamma \bH^\top\btR^{-1}=\\%
&\left( \bF^\top\btQ^{-1}\bF \right)^{-1} \bH^\top\btR^{-1}.
\end{align*}

Hence, \cref{eq:obspart} is given by
\begin{align}
&\Lambda \bF^\top\btQ^{-1}\bF %
\left( \bF^\top\btQ^{-1}\bF \right)^{-1}%
\bH^\top\btR^{-1}%
\left( \obs_k-\obsmod(\state_{k|k-1}) \right) = \nonumber\\
&\Lambda\bH^\top\btR^{-1}\left( \obs_k-\obsmod(\state_{k|k-1}) \right).\label{eq:obspartsimple}
\end{align}

Thus, reinserting \cref{term:C} and \cref{eq:obspartsimple} into \cref{eq:combinedsteps} gives
\begin{align}
\state_{k|k+1} =&\state_{k|k-1} + \Lambda\bigg[ \bH^\top\btR^{-1} \left( \obs_k-\obsmod(\state_{k|k-1}) \right) + \nonumber\\
& \bF^\top\btQ^{-1}(\state_{k+1|k+1} - \dynmod(\state_{k|k-1})) \bigg] = \nonumber\\
&\state_{k|k-1} + \Lambda\left[ \bH^\top(\bR+\bOmega_\obsmod)^{-1} \bigg( \obs_k-\obsmod(\state_{k|k-1}) \right) + \nonumber\\
& \bF^\top(\bQ+\bOmega_\dynmod)^{-1}(\state_{k+1|k+1} - \dynmod(\state_{k|k-1})) \bigg]\label{eq:finalcombined}.
\end{align}

Hence, with fixed $\bOmega_\obsmod,\bOmega_\dynmod, \bP_{k-1}$ and $\bP_k$, \cref{eq:finalcombined} is identical to the smoothing step in \cite{bellIteratedKalmanSmoother1994a} and thus minimizes the same loss, which, for one time step, is given by
\begin{align}
\mathcal{L} =~& (\state_{k-1}-\state_{k-1|k-1})^\top\bP_{k-1|k-1}^{-1}(\state_{k-1}-\state_{k-1|k-1})+\nonumber\\
&(\obs_k-\obsmod(\state_{k}))^\top(\bR+\bOmega_\obsmod)^{-1}(\obs_k-\obsmod(\state_{k})) + \nonumber\\
&(\state_{k}-\dynmod(\state_{k-1}))^\top(\bQ+\bOmega_\dynmod)^{-1}(\state_{k}-\dynmod(\state_{k-1})).
\end{align}
